\documentclass[12pt]{article}
\usepackage{amsfonts}
\usepackage{amsmath}
\usepackage{graphicx}

\input{tcilatex}
\begin{document}

\author{Ozlem Ersoy Hepson$^{1}$ and Idris Dag$^{2}$ \\
Department of Mathematics-Computer$^{1}$\\
Computer Engineering Department$^{2}$\\
Eski\c{s}ehir Osmangazi \"{U}niversity, Eski\c{s}ehir, Turkey,}
\date{}
\title{{\large \textbf{Numerical Solution of Singularly Perturbed Problems
via both Galerkin and Subdomain Galerkin methods}}}
\maketitle

\begin{abstract}
In this paper, numerical solutions of singularly perturbed boundary value
problems are given by using variants of finite element method. Both Galerkin
and subdomain Galerkin method based on quadratic B-spline functions are
applied over the geometrically graded. Results of some text problems are
compated with analytical solutions of the singularly perturbed problem

\end{abstract}

\section{Introduction}

This paper contains numerical solutions of one dimensional singularly
perturbation problems%
\begin{equation}
-\varepsilon u^{\prime \prime }+p(x)u^{\prime }+q(x)u=f(x),\text{ \ }x\in
\lbrack 0,1]  \label{1}
\end{equation}%
with boundary conditions%
\begin{equation}
u\left( 0\right) =\lambda \text{ and }u\left( 1\right) =\beta ,\text{ }%
\lambda ,\text{ }\beta \in 
\mathbb{R}
\label{2}
\end{equation}%
where $\varepsilon $ is a small positive parameter, $p(x),$ $q(x),$ $f(x)$
are sufficiently smooth functions with $p(x)\geq p^{\ast }>0,$ $q(x)\geq
q^{\ast }>0.$ These problems depend on $\varepsilon $ in such a way that the
solution varies rapidly in some parts and varies slowly in some other parts.
So, typically there are thin transition layers where the solutions can jump
abruptly, while away from the layers the solution behaves regularly and vary
slowly. The numerical treatment of the singular perturbation problems is far
from the trivial because of the boundary layer behavior of the solution.
There are a wide variety of asymptotic techniques for solving singular
perturbation problems.

These problems occur in many areas of engineering and applied mathematics
such as chemical reactor theory, optimal control, quantum mechanics, fluid
mechanics, reaction-diffusion process, aerodynamics, heat transport problems
with large Peclet numbers and Navier--Strokes flows withlarge Reynolds
numbers etc.

Many authors have studied on this problem and tried to overcome the
above-mentioned difficulties. M.K. Kadalbajoo and Vikas Gupta \cite%
{kadalbajoo} proposed B-spline collocation method on a non-uniform mesh of
Shishkin type.to solve singularly perturbed two-point boundary value
problems with a turning point exhibiting twin boundary layers. J.Vigo-Aguiar
and S.Natesan \cite{aguiar} consider a class of singularly perturbed
two-point boundary-value problems for second-order ordinary differential
equations. They suggested an iterative non-overlapping domain decomposition
method in order to obtain numerical solution to these problems. Tirmizi et
al. \cite{Tirmizi} have proposed a generalized scheme based on quartic
non-polynomial spline functions in order to designed for numerical solution
of singularly perturbed two-point boundary-value problems. D.J.Fyfe \cite%
{fyfe} used cubic splines on equal and unequal intervals and compared the
results. He observed that very little advantage is gained by using unequal
intervals. M.K.Kadalbajoo and K.C.Patidar \cite{kadalbajoo1} gave some
difference schemes using spline in tension. They showed that these methods
are second-order accurate. Employing coordinate stretching a
Galerkin-spectral method is applied to the singularly perturbed boundary
value problems by W.Liu and T.Tang \cite{wenbin}. G.Beckett and
J.A.Mackenzie \cite{beckett} gave a $p$th order Galerkin finite element
method on a non-uniform grid. In their study the grid is constructed by
equidistributing a strictly positive monitor function. After the appropriate
selection of the monitor function parameters they obtained insensitive
numerical solution.

The definitions of B-splines over the geomertically graded mesh was given in
reference \cite{dag}. Dag and Sahin \cite{Sahin} have set up the finite
element method employing the quadratic and the cubic B-splines to form the
trial function. In this article, we used the finite element method with the
quadratic B-splines. After giving the expressions of the mentioned B-splines
over the geometrically graded mesh we applied the quadratik Galerkin and
quadratik subdomain Galerkin method to Eq.(\ref{1}).

Briefly, outline is as follows. In Section 2, numerical methods are given.
Numerical experiments are carried out for one test problem and errors of
those methods are compared in Section 3. Finally conclusion is given in last
section.

\section{B-spline Galerkin Methods}

For numerical purpose, let us divide the solution domain $\left[ 0,1\right] $
into subintervals by the knots $x_{m}$ such that%
\begin{equation*}
0=x_{0}<x_{1}<\cdot \cdot \cdot <x_{N}=1
\end{equation*}%
where $x_{m+1}=x_{m}+h_{m}$ and $h_{m}$ is the size of interval $\left[
x_{m},x_{m+1}\right] $ with relation $h_{m}=\sigma h_{m-1}.$ Here $\sigma $
is mesh ratio constant.

To construct the geometrically graded mesh, determination of the first
element size $h_{0}$ is necessary. Since%
\begin{equation*}
h_{0}+h_{1}+\cdot \cdot \cdot +h_{N-1}=1
\end{equation*}%
it is easy to write 
\begin{equation*}
h_{0}=\frac{1}{1+\sigma +\sigma ^{2}+\cdot \cdot \cdot +\sigma ^{N-1}}.
\end{equation*}

This partition will be uniform if the mesh ratio $\sigma $ is taken as
unity. To obtain finer mesh at the left boundary, $\sigma $ must be chosen
as $\sigma >1.$ On the other hand, to make the mesh size smaller at the
right boundary, $\sigma $ must be chosen as $\sigma <1$ . Mentioned
selection of $\sigma $ will be done by experimentally.

\subsection{Quadratic B-spline Galerkin method (QM)}

The expression of the quadratic B-splines over the geometrically graded mesh
may be given in the following form \cite{dag}:%
\begin{equation}
\begin{tabular}{l}
$Q_{m-1}$ \\ 
$Q_{m}$ \\ 
$Q_{m+1}$%
\end{tabular}%
=\dfrac{1}{h_{m}^{2}}\left \{ \ 
\begin{array}{l}
(h_{m}-\xi )^{2}\sigma , \\ 
h_{m}^{2}+2h_{m}\sigma \xi -\left( 1+\sigma \right) \xi ^{2}, \\ 
\xi ^{2}%
\end{array}%
\  \right.  \label{3}
\end{equation}

where $\xi =x-x_{m}$ and $0\leq \xi \leq h_{m}.$ A quadratic B-spline covers 
$3$ elements. Any quadratic B-spline $Q_{m}$ and its derivatives vanish
outside of the interval $\left[ x_{m-1},x_{m+2}\right] $ and therefore an
element is covered by $3$ successive quadratic B-splines. The set of the
quadratic B-splines $\left \{ Q_{-1},Q_{0},...,Q_{N}\right \} $ forms a
basis for the functions defined on the solution domain \cite{prenter}.
Thence, an approximation $u_{N}$ to the analytical solution $u$ can be
written as%
\begin{equation}
u_{N}=\sum \limits_{m=-1}^{N}\delta _{m}Q_{m}  \label{4}
\end{equation}%
where $\delta _{m}$ are unknown parameters. By the substitution of the value
of $Q_{m}$ at the knots $x_{m}$ in Eq.(\ref{4}), the nodal value $u$\ and
its derivative $u^{\prime }$ are expressed in terms of $\delta _{m}$ by\ 
\begin{equation}
\begin{array}{ll}
u_{m}=u(x_{m}) & =\sigma \delta _{m-1}+\delta _{m},\vspace{0.2cm} \\ 
u_{m}^{\prime }=u^{\prime }(x_{m}) & =\dfrac{2\sigma }{h_{m}}(\delta
_{m}-\delta _{m-1}).%
\end{array}
\label{5}
\end{equation}%
Both sides of the weight function by multiplying the differential equation
and is integrated over the range $[x_{m},x_{m+1}]$ following equation is
obtained.%
\begin{equation*}
-\varepsilon vu^{\prime \prime }(x)+p(x)vu^{\prime }(x)+q(x)vu(x)=f(x)
\end{equation*}%
partial integration is applied to the first term in the above integral is
obtained as follows.%
\begin{equation*}
\dint \limits_{x_{m}}^{x_{m+1}}(-\varepsilon v^{\prime }u^{\prime
}(x)+vp(x)u^{\prime }+vq(x)u(x))dx-\varepsilon vu^{\prime }(x)\overset{%
x_{m+1}}{\underset{x_{m}}{\mid }}-\dint \limits_{x_{m}}^{x_{m+1}}vf(x)dx=0
\end{equation*}%
weight functions selected as $Q_{j},$ $j=m-1,m,m+1$ and used (\ref{4}),
following integral is obtained.%
\begin{equation}
\dsum \limits_{j=m-1}^{m+1}[-\varepsilon \dint \limits_{0}^{h_{m}}(\phi
_{i}^{\prime }\phi _{j}^{\prime }+p\phi _{i}\phi _{j}^{\prime }+q\phi
_{i}\phi _{j})d\xi ]\delta _{j}-\varepsilon \phi _{i}\phi _{j}^{\prime }%
\overset{h_{m}}{\underset{0}{\mid }}\delta _{j}-\dint
\limits_{0}^{h_{m}}\phi _{i}f(x_{m}+\xi )d\xi =0  \label{3.6}
\end{equation}%
So, following values can be computed.%
\begin{equation*}
\begin{array}{lll}
a_{ij}=\dint \limits_{0}^{h_{m}}\phi _{i}^{\prime }\phi _{j}^{\prime }d\xi ,
&  & r_{ij}=\phi _{i}\phi _{j}^{\prime }\overset{h_{m}}{\underset{0}{\mid }},
\\ 
&  &  \\ 
b_{ij}=\dint \limits_{0}^{h_{m}}\phi _{i}\phi _{j}^{\prime }d\xi , &  & 
f_{i}=\dint \limits_{0}^{h_{m}}\phi _{i}f(x_{m}+\xi )d\xi , \\ 
\multicolumn{1}{c}{} & \multicolumn{1}{c}{} & \multicolumn{1}{c}{} \\ 
\multicolumn{1}{c}{c_{ij}=\dint \limits_{0}^{h_{m}}\phi _{i}\phi _{j}d\xi ,}
& \multicolumn{1}{c}{} & \multicolumn{1}{c}{}%
\end{array}%
\end{equation*}%
Where $i,j=m-1,m,m+1$ and 
\begin{eqnarray*}
A^{(m)} &=&\frac{2}{3h_{m}}%
\begin{bmatrix}
2\alpha ^{2} & \alpha (1-\alpha ) & -\alpha \\ 
\alpha (1-\alpha ) & 2(1-\alpha +\alpha ^{2}) & \alpha -2 \\ 
-\alpha & \alpha -2 & 2%
\end{bmatrix}
\\
B^{(m)} &=&%
\begin{bmatrix}
\frac{-\alpha ^{2}}{2} & \frac{1}{6}\alpha \left( 3\alpha -1\right) & \frac{%
\alpha }{6} \\ 
-\frac{1}{6}\alpha \left( 3\alpha +5\right) & \frac{1}{2}\alpha ^{2}-\frac{1%
}{2} & \frac{5}{6}\alpha +\frac{1}{2} \\ 
\frac{-\alpha }{6} & \frac{1}{6}\alpha -\frac{1}{2} & \frac{1}{2}%
\end{bmatrix}
\\
C^{(m)} &=&h_{m}%
\begin{bmatrix}
\frac{1}{5}\alpha ^{2}h_{m} & \frac{1}{30}\alpha h_{m}\left( 4\alpha
+9\right) & \frac{\alpha }{30} \\ 
\frac{1}{30}\alpha h_{m}\left( 4\alpha +9\right) & \frac{8}{15}\alpha ^{2}+%
\frac{11}{5}\alpha +\frac{8}{15} & \frac{3}{10}\alpha +\frac{2}{15} \\ 
\frac{\alpha }{30} & \frac{3}{10}\alpha +\frac{2}{15} & \frac{1}{5}h_{m}%
\end{bmatrix}
\\
R^{(m)} &=&\frac{1}{h_{m}}%
\begin{bmatrix}
2\alpha ^{2} & 0 & 0 \\ 
2\alpha (3\alpha +2) & -4\alpha & 2\alpha \\ 
0 & 0 & 0%
\end{bmatrix}%
\end{eqnarray*}%
and%
\begin{equation*}
F^{(m)}=\frac{1}{h_{m}^{2}}%
\begin{bmatrix}
\vartheta _{1} & \vartheta _{1} & \vartheta _{1} \\ 
\vartheta _{2} & \vartheta _{2} & \vartheta _{2} \\ 
\vartheta _{3} & \vartheta _{3} & \vartheta _{3}%
\end{bmatrix}%
\end{equation*}%
Where%
\begin{eqnarray*}
\vartheta _{1} &=&-\alpha \left( 2h_{m}-2e^{h_{m}}+h_{m}^{2}+2\right) \\
\vartheta _{2} &=&2(\alpha +1)(1-e^{h_{m}})+2h_{m}(\alpha
+e^{h_{m}})-h_{m}^{2}(1-\alpha e^{h_{m}}) \\
\vartheta _{3} &=&\left( e^{h_{m}}\left( h_{m}^{2}-2h_{m}+2\right) -2\right)
\end{eqnarray*}%
Defined in terms of local matrices $A^{(i)},$ $B^{(i)},$ $C^{(i)},$ $R^{(i)}$
and $F^{(i)}$, equation (\ref{3.6}) can be represented following form.%
\begin{equation*}
(-\varepsilon A^{(i)}+pB^{(i)}+qC^{(i)}-\varepsilon R^{(i)})\delta
^{(i)}=F^{(i)}
\end{equation*}%
Where%
\begin{equation*}
\delta ^{(i)}=(\delta _{m-1}^{(i)},\delta _{m}^{(i)},\delta
_{m+1}^{(i)},\delta _{mi+2}^{()}),\text{ }%
F^{(i)}=(f_{m-1}^{(i)},f_{m}^{(i)},f_{m+1}^{(i)},f_{m+2}^{(i)})
\end{equation*}%
By combining the local matrices which is defined on $[x_{i},x_{i+1}],$ $%
i=0,\ldots ,N-1$, the global system in the range of $[x_{0},x_{N}]$ can be
defined as follows.%
\begin{equation}
(-\varepsilon A+pB+qC-\varepsilon R)\delta =F  \label{3.7}
\end{equation}%
The matrix of $A$ is%
\begin{equation*}
A=\left[ 
\begin{array}{cccccccccc}
\sigma _{0,0}^{(0)} & \sigma _{0,1}^{(0)} & \sigma _{0,2}^{(0)} &  &  &  & 
&  &  &  \\ 
\sigma _{1,0}^{(0)} & \sigma _{1,1}^{\ast (1)} & \sigma _{1,2}^{\ast (1)} & 
\sigma _{1,3}^{\ast (1)} &  &  &  &  &  &  \\ 
\sigma _{2,0}^{(0)} & \sigma _{2,1}^{\ast (1)} & \sigma _{2,2}^{\ast (2)} & 
\sigma _{2,3}^{\ast (2)} & \sigma _{2,4}^{\ast (2)} &  &  &  &  &  \\ 
& \sigma _{3,0}^{\ast (1)} & \sigma _{3,1}^{\ast (2)} & \sigma _{3,2}^{\ast
(3)} & \sigma _{3,3}^{\ast (3)} & \sigma _{3,4}^{\ast (3)} &  &  &  &  \\ 
&  & \ddots & \ddots & \ddots & \ddots & \ddots &  &  &  \\ 
&  &  & \sigma _{i,i-2}^{\ast (i-1)} & \sigma _{i,i-1}^{\ast (i)} & \sigma
_{i,i}^{\ast (i+1)} & \sigma _{i,i+1}^{\ast (i+1)} & \sigma _{i,i+2}^{\ast
(i+1)} &  &  \\ 
&  &  &  & \ddots & \ddots & \ddots & \ddots & \ddots &  \\ 
&  &  &  &  & \sigma _{n-1,n-3}^{\ast (n-2)} & \sigma _{n-1,n-2}^{\ast (n-1)}
& \sigma _{n-1,n-1}^{\ast (n)} & \sigma _{n-1,n}^{\ast (n)} & \sigma
_{n-1,n+1}^{(n)} \\ 
&  &  &  &  &  & \sigma _{n,n-2}^{\ast (n-1)} & \sigma _{n,n-1}^{\ast (n-1)}
& \sigma _{n,n}^{\ast (n-1)} & \sigma _{n,n+1}^{(n)} \\ 
&  &  &  &  &  &  & \sigma _{n+1,n-1}^{(n)} & \sigma _{n+1,n}^{(n)} & 
a_{n+1,n+1}^{(n)}%
\end{array}%
\right]
\end{equation*}%
Where%
\begin{eqnarray*}
&&%
\begin{array}{lll}
\sigma _{1,1}^{\ast (1)}=\sigma _{1,1}^{(0)}+\sigma _{1,1}^{(1)}, &  & 
\sigma _{1,2}^{\ast (1)}=\sigma _{1,2}^{(0)}+\sigma _{1,2}^{(1)}, \\ 
&  &  \\ 
\sigma _{2,1}^{\ast (1)}=\sigma _{2,1}^{(0)}+\sigma _{2,1}^{(1)}, &  & 
\sigma _{2,2}^{\ast (2)}=\sigma _{2,2}^{(0)}+\sigma _{2,2}^{(1)}+\sigma
_{2,2}^{(2)}, \\ 
&  &  \\ 
\sigma _{2,3}^{\ast (2)}=\sigma _{2,3}^{(2)}+\sigma _{2,3}^{(1)}, &  & 
\sigma _{i,i-1}^{\ast (i)}=\sigma _{i,i-1}^{(i-1)}+\sigma _{i,i-1}^{(i)},%
\end{array}
\\
&&%
\begin{array}{lll}
\sigma _{i,i}^{\ast (i)}=\sigma _{i,i}^{(i-1)}+\sigma _{i,i}^{(i)}+\sigma
_{i,i}^{(i+1)}, &  & \sigma _{i,i+1}^{\ast (i)}=\sigma _{i,i+1}^{(i)}+\sigma
_{i,i+1}^{(i+1)}, \\ 
&  &  \\ 
\sigma _{n-1,n-2}^{\ast (n-1)}=\sigma _{n-1,n-2}^{(n-2)}+\sigma
_{n-1,n-2}^{(n-1)}, &  & \sigma _{n-1,n-1}^{\ast (n)}=\sigma
_{n-1,n-1}^{(n-2)}+\sigma _{n-1,n-1}^{(n-1)}+\sigma _{n-1,n-1}^{(n)}, \\ 
&  &  \\ 
\sigma _{n-1,n}^{\ast (n)}=\sigma _{n-1,n}^{(n-1)}+\sigma _{n-1,n}^{(n)}, & 
& \sigma _{n,n-1}^{\ast (n-1)}=\sigma _{n,n-1}^{(n-1)}+\sigma _{n,n-1}^{(n)},
\\ 
&  &  \\ 
\sigma _{n,n}^{\ast (n-1)}=\sigma _{n,n}^{(n-1)}+\sigma _{n,n}^{(n)}. &  & 
\end{array}%
\end{eqnarray*}%
$B,$ $C,$ $R$ matrices are obtained similarly. Also the matrix of $F$ is
computed as follows.%
\begin{equation*}
F=\left[ 
\begin{array}{l}
f_{0}^{0} \\ 
f_{1}^{0}+f_{1}^{1} \\ 
f_{2}^{0}+f_{2}^{1}+f_{2}^{2} \\ 
\vdots \\ 
f_{i}^{i-1}+f_{i}^{i}+f_{i}^{i+1} \\ 
\vdots \\ 
f_{n-2}^{n-1}+f_{n-2}^{n-2}+f_{n-2}^{n-3} \\ 
f_{n-1}^{n-1}+f_{n-1}^{n-2} \\ 
f_{n}^{n-1}%
\end{array}%
\right]
\end{equation*}%
The matrix system (\ref{3.7}) has $N+1$ equations and $N+3$ unknowns. In
order to solve this system, the numbers of equations and unknowns must be
equal. From the boundary conditions (\ref{2}) and Eq.(\ref{5}) it is easy to
write%
\begin{equation*}
\delta _{-1}=\frac{\lambda -\delta _{0}}{\alpha },\text{ \ }\delta
_{N}=\beta -\alpha \delta _{N-1}.
\end{equation*}%
Using these equalities, $\delta _{-1}$ and $\delta _{N}$ can be eliminated
from the system and then matrix equation (\ref{7.1}) can be solved with
Thomas algorithm. Substituting the obtained parameters $\delta _{m}$ in Eq.(%
\ref{5}), the numerical solution is found at the knots $x_{m}.$

\subsection{Quadratic B-spline Subdomain Galerkin method (QM)}

If on each side of the equation (\ref{1}), multiplied by weight function $%
V_{n}$ and is integrated over the range $[x_{m},x_{m+1}]$ then following
integral obtained.%
\begin{equation*}
\dint \limits_{x_{0}}^{x_{n}}[-\varepsilon u^{\prime \prime }(x)+pu^{\prime
}(x)+qu(x)-f(x)]dx=0
\end{equation*}%
Here the wight function is%
\begin{equation*}
V_{n}=\left \{ 
\begin{array}{cc}
1, & x_{m}\leq x<x_{m+1} \\ 
0, & \text{other case}%
\end{array}%
\right.
\end{equation*}%
partial integration is applied to the first term in the above integral is
obtained as follows.%
\begin{equation*}
-\varepsilon u^{\prime }(x)\overset{x_{m+1}}{\underset{x_{m}}{\mid }}+pu(x)%
\overset{x_{m+1}}{\underset{x_{m}}{\mid }}+q\dint%
\limits_{x_{m}}^{x_{m+1}}u(x)dx=\dint \limits_{x_{m}}^{x_{m+1}}f(x)dx
\end{equation*}%
By the substitution of the nodal value $u$ and its derivative $u^{\prime }$
in last equation is expressed following integral.%
\begin{equation}
\begin{array}{l}
\lbrack -\varepsilon (\dsum \limits_{J=m-1}^{m+1}\phi _{i}^{\prime }\overset{%
hm}{\underset{0}{\mid }})\delta _{j})+p(x_{m}+\xi )(\dsum
\limits_{J=m-1}^{m+1}\phi _{i}\overset{hm}{\underset{0}{\mid }}\delta
_{j})+q(x_{m}+\xi )\dsum \limits_{J=m-1}^{m+1}\dint
\limits_{x_{m}}^{x_{m+1}}\phi _{i}\delta _{j}dx \\ 
=\dint \limits_{x_{m}}^{x_{m+1}}f(x_{m}+\xi )dx%
\end{array}
\label{4.1}
\end{equation}%
With the help of the division points values, Qauadratik B-spline shape
functions defined on geometrically increasing intervals $[x_{m},x_{m+1}]$%
\begin{eqnarray}
\dint \limits_{x_{m}}^{x_{m+1}}u^{\prime \prime }(x)dx &=&u^{\prime }(x)%
\overset{x_{m+1}}{\underset{x_{m}}{\mid }}=\frac{2\alpha }{h_{m+1}}\delta
_{m+1}+(-\frac{2\alpha }{h_{m+1}}-\frac{2\alpha }{h_{m}})\delta _{m}+\frac{%
2\alpha }{h_{m}}\delta _{m-1}  \label{5.1} \\
\dint \limits_{x_{m}}^{x_{m+1}}u^{\prime }(x)dx &=&u(x)\overset{x_{m+1}}{%
\underset{x_{m}}{\mid }}=\delta _{m+1}+(\alpha -1)\delta _{m}-\alpha \delta
_{m-1}  \label{5.2} \\
\dint \limits_{x_{m}}^{x_{m+1}}u(x)dx &=&\dint \limits_{0}^{h_{m}}(\dsum
\limits_{J=m-1}^{m+1}\phi _{j}\delta _{j})d\xi =\delta _{m-1}\dint
\limits_{x_{m}}^{x_{m+1}}\phi _{m-1}dx+\delta _{m}\dint
\limits_{x_{m}}^{x_{m+1}}\phi _{m}dx+\delta _{m+1}\dint
\limits_{x_{m}}^{x_{m+1}}\phi _{m+1}dx  \label{5.3}
\end{eqnarray}%
The value of $Q_{m-1,}Q_{m},$\bigskip $Q_{m+1}$ substitute in (\ref{5.3}) is
computed as follows.%
\begin{eqnarray*}
\dint \limits_{x_{m}}^{x_{m+1}}\phi _{m-1}dx &=&\frac{1}{3}\alpha h_{m} \\
\dint \limits_{x_{m}}^{x_{m+1}}\phi _{m}dx &=&\frac{2}{3}h_{m}\left( \alpha
+1\right) \\
\dint \limits_{x_{m}}^{x_{m+1}}\phi _{m+1}dx &=&\frac{1}{3}h_{m}
\end{eqnarray*}%
If we replace integrals calculated by (\ref{5.3}) the following results are
obtained.%
\begin{equation*}
\dint \limits_{x_{m}}^{x_{m+1}}u(x)dx=\frac{1}{3}\alpha h_{m}\delta _{m-1}+%
\frac{2}{3}h_{m}\left( \alpha +1\right) \delta _{m}+\frac{1}{3}h_{m}\delta
_{m+1}
\end{equation*}

When we apply the Galerkin method to the system (\ref{4.1}), (\ref{5.1}), (%
\ref{5.2}) and (\ref{5.3}) derivatives are replaced with their equals
obtained from Eq.(\ref{4.1}). This substitution yields the following system:%
\begin{equation*}
\begin{array}{l}
(-\dfrac{2\alpha \varepsilon }{h_{m}}+\alpha (\dfrac{1}{3}%
h_{m}q(x)-p(x)))\delta _{m-1}+(\dfrac{2\varepsilon }{h_{m}}(1+\alpha
)+(\alpha -1)p(x)+\dfrac{2}{3}h_{m}\left( \alpha +1\right) q(x))\delta _{m}
\\ 
\\ 
+(-\dfrac{2\varepsilon }{h_{m}}+p(x)+\dfrac{1}{3}h_{m}q(x))\delta
_{m+1}=\dint \limits_{x_{m}}^{x_{m+1}}f(x)dx%
\end{array}%
\end{equation*}%
With necessary operations, this system can be written in matrix form as%
\begin{equation}
\mathbf{AX=F}  \label{7}
\end{equation}%
where%
\begin{equation}
A=\left[ 
\begin{array}{cccccccc}
\alpha _{01} & \alpha _{02} & \alpha _{03} &  &  &  &  &  \\ 
& \alpha _{11} & \alpha _{12} & \alpha _{13} &  &  &  &  \\ 
&  & \alpha _{21} & \alpha _{22} & \alpha _{23} &  &  &  \\ 
&  &  & \alpha _{31} & \alpha _{32} & \alpha _{32} &  &  \\ 
&  &  & \ddots & \ddots & \ddots &  &  \\ 
&  &  &  &  & \alpha _{n1} & \alpha _{n2} & \alpha _{n3}%
\end{array}%
\right] ,  \label{7.1}
\end{equation}%
where%
\begin{equation*}
\begin{array}{l}
\alpha _{m1}=-\dfrac{2\alpha \varepsilon }{h_{m}}-\alpha p_{m}+\dfrac{1}{3}%
\alpha h_{m}q_{m}, \\ 
\\ 
\alpha _{m2}=\dfrac{2\varepsilon }{h_{m}}(1+\alpha )+(\alpha -1)p_{m}+\dfrac{%
2}{3}h_{m}\left( \alpha +1\right) q_{m} \\ 
\\ 
\alpha _{m3}=-\dfrac{2\alpha \varepsilon }{\alpha h_{m}}+p_{m}+\dfrac{1}{3}%
h_{m}q_{m}%
\end{array}%
\end{equation*}%
and%
\begin{eqnarray*}
X &=&\left[ 
\begin{array}{cccccc}
\delta _{-1}, & \delta _{0}, & \delta _{1}, & \cdots , & \delta _{n-1}, & 
\delta _{n}%
\end{array}%
\right] ^{T} \\
F &=&\left[ 
\begin{array}{cccccc}
f_{0}, & f_{1}, & f_{2}, & \cdots , & f_{n-1}, & f_{n}%
\end{array}%
\right] ^{T}
\end{eqnarray*}%
\begin{equation*}
\text{ }f_{m}=f(x_{m}),\text{ }m=0,1,\cdots ,N
\end{equation*}

The matrix system (\ref{7}) has $N+1$ equations and $N+3$ unknowns. In order
to solve this system, the numbers of equations and unknowns must be equal.
From the boundary conditions (\ref{2}) and Eq.(\ref{5}) it is easy to write%
\begin{equation*}
\delta _{-1}=\frac{\lambda -\delta _{0}}{\alpha },\text{ \ }\delta
_{N}=\beta -\alpha \delta _{N-1}.
\end{equation*}%
Using these equalities, $\delta _{-1}$ and $\delta _{N}$ can be eliminated
from the system and then matrix equation (\ref{7.1}) can be solved with
Thomas algorithm. Substituting the obtained parameters $\delta _{m}$ in Eq.(%
\ref{5}), the numerical solution is found at the knots $x_{m}.$

\section{Numerical Experiments}

We have tested the accuracy of the numerical methods on two examples. Errors
are measured with the norm 
\begin{equation*}
L_{\infty }=\left \vert u-u_{N}\right \vert _{\infty }=\max_{j}\left \vert
u_{j}-\left( u_{N}\right) _{j}\right \vert .
\end{equation*}%
Since the boundary layers are at the right boundary in both examples, in
order to minimize the error, we have searched the interval $\left(
0,1\right) $ for the best choice of \ the mesh ratio $\sigma .$ Solution
profiles are illustrated in Figs. 1-4 for the first example and in Figs. 5-8
for the second example. These figures are graphed for $N=20$ and two
different $\varepsilon .$ In order to see the success of the numerical
methods more clear, exact solutions and obtained results are illustrated
together in all figures. In all figures, continuous line is used for the
exact solutions and the lines$\cdots $o$\cdots $o$\cdots $, $\cdots $+$%
\cdots $+$\cdots $ are used for QM and CM respectively. Using uniform mesh
leads to oscillations, seen in Figs. 1, 3, 5 and 7, in solution profiles
because of the boundary layer. As observed from Figs. 2, 4, 6 and 8, after
the best choice of mesh ratio $\sigma ,$ these oscillations disappear.\
Using various $\varepsilon $ and $N,$ calculated numerical errors are
tabulated and compared in Table 1 and Table 2 for the first and the second
examples respectively.

\textbf{Example }Our example is 
\begin{equation*}
\begin{tabular}{c}
$-\varepsilon u^{\prime \prime }+u^{\prime }=\exp (x),$ \\ 
\\ 
$u(0)=u(1)=0$%
\end{tabular}%
\end{equation*}%
with the exact solution%
\begin{equation*}
u(x)=\frac{1}{1-\varepsilon }\left[ \exp (x)-\frac{1-\exp (1-1/\varepsilon
)+\left( \exp (1)-1\right) \exp \left( \left( x-1\right) /\varepsilon
\right) }{1-\exp (-1/\varepsilon )}\right] .
\end{equation*}%
taken from \cite{lorenz}.

\newpage

\bigskip

\section{Conclusion}

Quadratic and cubic B-spline algorithms are applied to singularly perturbed
problems. Difficulties arised from the modelling of the boundary layers in
numerical methods are tried to overcome by using B-splines over the
geometrically graded mesh. Simplicity of the adaptation of B-splines\ and
obtaining acceptable solutions can be noted as advantages of given numerical
methods. Consequently, in getting the numerical solution of the differential
equations having boundary layers, B-spline collocation methods over the
geometrically graded mesh are advisable.

\end{document}